\documentclass[sn-mathphys]{sn-jnl}

\definecolor{mon}{rgb}{10,126,140}
\definecolor{royal}{rgb}{1, 0, 0}
\definecolor{blue}{rgb}{0, 0, 1}
\definecolor{raspberryrose}{rgb}{0.7, 0.27, 0.42}
\definecolor{olive}{rgb}{0.5, 0.5, 0.0}
\definecolor{blush}{rgb}{0.87, 0.36, 0.51}
 \definecolor{persianblue}{rgb}{0.11, 0.22, 0.73}
\definecolor{junebud}{rgb}{0.74, 0.85, 0.34}
\definecolor{heartgold}{rgb}{0.5, 0.5, 0.0}

\def\1{\mathds{1}}

\jyear{2021}%

\theoremstyle{thmstyleone}%
\newtheorem{theorem}{Theorem}
\newtheorem{proposition}[theorem]{Proposition}%
\newtheorem{lemma}[theorem]{Lemma}%
\newtheorem{corollary}[theorem]{Corollary}%

\theoremstyle{thmstyletwo}%
\newtheorem{problem}{Problem}%

\theoremstyle{thmstylethree}%

\begin{document}

\title[
Convexity properties related to Gauss
hypergeometric function]{
Convexity properties related to Gauss
hypergeometric function}


\author*[1,2]{\fnm{Mohamed} \sur{Bouali
}}\email{bouali25@laposte.net}



\affil*[1]{\orgdiv{Department of mathematics}, \orgname{Preparatory Institute for Engineering Studies of Tunis}, \orgaddress{\street{Jawaher Lal Naherou}, \city{Montfleury}, \postcode{1089}, \state{Tunis}, \country{Tunisia}}}

\affil[2]{\orgdiv{Department of mathematics}, \orgname{Faculty of Sciences of Tunis}, \orgaddress{\street{Rue Tol\`ede}, \city{El-Manar}, \postcode{2092}, \state{Tunis}, \country{Tunisia}}}



\abstract
{We investigate the convexity property on $(0,1)$ of the functions  $\varphi_{a,b,c}$ and $1/\varphi_{a,b,c}$, where $$\varphi_{a,b,c}(x)= \frac{c-\log(1-x)}{\,_2F_1(a,b,a+b,x)},$$
whenever $a,b\geq 0$ and $a+b\leq 1$.

We Show that $\varphi_{a,b,c}$ (respectively $1/\varphi_{a,b,c}$) is strictly convex on $(0,1)$ if and only if $c\leq -2\gamma-\psi(a)-\psi(b),$ (respectively $c\geq\alpha_0$) and $\varphi_{a,b,c}$ (respectively $1/\varphi_{a,b,c}$) is strictly concave on $(0,1)$ if and only if $c\geq c(a,b)$ (respectively $c\in[\delta_-,\delta_+]$), where $\psi$ is the Polygamma function. This generalizes some problems posed by Yang and Tian and complete the study of convexity properties of functions studied by the author in [bouali]. As applications of the convexity and concavity, we establish among other inequalities, that for all $x\in(0,1)$, $a,b\in[0,1]$, $a+b\leq 1$ and $c\geq c(a,b)$
$$c+\frac{\Gamma(a)\Gamma(b)}{\Gamma(a+b)}\leq \frac{c-\log(1-x)}{\,_2F_1(a,b,a+b,x)}+\frac{c-\log(x)}{\,_2F_1(a,b,a+b,1-x)}\leq\frac{(2c+2\log 2)}{\,_2{F}_1(a,b;a+b;1/2)},$$
and
for all $x\in(0,1)$, $a,b\in[0,1]$, $a+b\leq 1$ and $c\in [\delta_-,\delta_+]$
$$\frac1c+\frac{\Gamma(a+b)}{\Gamma(a)\Gamma(b)}\leq \frac{\,_2F_1(a,b,a+b,x)}{c-\log(1-x)}+\frac{\,_2F_1(a,b,a+b,1-x)}{c-\log(x)}\leq\frac{\,_2{F}_1(a,b;a+b;1/2)}{(2c+2\log 2)}.$$.}

\keywords{  Convexity; Gauss hypergeometric function; inequalities.}


\pacs[MSC Classification]{26D07; 33C05; 33E05}

\maketitle

\bigskip

\section{Introduction}\label{sec1}
For real numbers $x, y > 0$, the gamma, beta and psi functions are defined as
$$\Gamma(x)=\int_0^\infty t^{x-1}e^{-t}dt,\quad B(x,y)=\frac{\Gamma(x)\Gamma(y)}{\Gamma(x+y)},\quad\psi(x)=\frac{\Gamma'(x)}{\Gamma(x)}.$$
The Gauss hypergeometric function is defined on $(-1,1)$ by the series
$${}_2F_1(a,b,c,x)=\sum_{n=0}^\infty\frac{(a)_n(b)_n}{(c)_n}\frac{x^n}{n!},$$
where $(a)_n=a(a+1)\cdots (a+n-1)=\Gamma(a+n)/\Gamma(a)$ is the Pochhammer symbol.
Let $\gamma=0.5772156649\cdots=-\psi(1)$ be the Euler-Mascheroni constant, and for $a\in(0,1)$, let $R(a)$ be defined by
$$R(a,b)=-2\gamma-\psi(a)-\psi(b),$$
which is called the Ramanujan constant in literature although it is actually a function of $a$ and $b$.
It is well known that $R(a,1-a)$ is essential in some fields of mathematics such as the zero-balanced Gaussian
hypergeometric functions ${}_2F_1(a,1-a,1,x)$, the theories of Ramanujan's modular equations and quasiconformal
mappings.

Specials functions and specially elliptic functions arises in numerous branches of mathematics such as geometric function theory
and quasi-conformal mappings, also in in physics, theory of mean values, number theory and other related fields, see for instance
\cite{byr, ber, bal, bha, ande, qiu, wang, alz}. Many authors were interested in studying monotonicity properties of functions related to elliptic functions ${\cal  K}(x)=(\pi/2) {}_2F_1(1/2,1/2,1,x)$ and ${\cal E}(x)=(\pi/2)  {}_2F_1(1/2,-1/2,1,x)$. Anderson et al. in \cite{ander} showed that the function defined on $(0,1)$ by
$$U_a(x)=\frac{a-\log(1-x^2)}{{\cal K}(x)},$$
is strictly decreasing if and only if $ a\geq\log 16$, and strictly increasing if and only if $a\leq 4$.
 We generalize this result in the following way. Let $a,b\geq 0$, then the function $\varphi_{a,b,c}$ is strictly decreasing on $(0,1)$ if and only if $c\geq 1/a+1/b$ and is strictly increasing on $(0,1)$ if and only if $c\leq R(a,b)$.

In a recent paper \cite{yan}, Yang and Tian  studied the closely related function
$$V_b(x)=\frac{{\cal K}(\sqrt x)}{b-(1/2)\log(1-x)},\quad c\geq 0.$$
They proved that $V_b$ is strictly concave on $(0,1)$ if and only if $b=4/3$. They also posed the following problem.

{\it Determine the best parameters $a$ and $b$ such that $V_a$ is convex on $(0, 1)$ and $1/V_b$ is concave on $(0, 1)$.}

In the recent paper \cite{alz1}, Alzer and Richards give an answer to the second problem.

Later in \cite{rich}, Richards and Smith extended the second problem to the generalized elliptic integral ${\cal K}_p, (p\geq 1)$. However, the second problem is studied by the author in [Bouali].

The objective of this paper is to extend these problems by studying convexity properties of some functions related to Gauss hypergeometric function.

\section{Main results}

\begin{theorem}\label{t2} For $a,b\geq 0$, $c\geq 0$, let
$$\varphi_{a,b,c}(x)=\frac{c-\log(1-x)}{\,_2F_1(a,b,a+b,x)}.$$
If $a+b\leq 1$, then
the function $\varphi_{a,b,c}$ is strictly convex on $(0,1)$ if and only if $c\leq R(a,b)$ and is strictly concave on $(0,1)$ if and only if
$c\geq c(a,b)$, where
$$c(a,b)=\frac{(a+b)(a+b-2 a b)(a+b+1)}{a b ((a+b+1)(a+b-2 a b)+ab(a+b))}.$$
\end{theorem}
For $a=b=1/2$, we retrieve the result of Alzer with $c(1/2,1/2)=16/5$ and the case $a+b=1$, $a,b\geq 0$ corresponds to the result of Kendall with
$c(a,1-a)=2(1-2a+2a^2)/a((1-a)(2-3a+3a^2))$.

To state the second main result, let us introduce some auxiliary functions.
For $a,b\geq 0$, we set
$$\varphi_\pm(x)=\log(1-x)+\omega_{\pm}(x),$$ where,
$$\omega_\pm(x)=\frac{-g(x)\pm\sqrt{\Delta(x)}}{2h(x)}.$$
$$h(x)=\frac{(ab)^2}{(a+b)(a+b+1)}(1-x)\,_2F_1(a+1,b+1,a+b+2,x)+\frac{ab}{a+b}\,_2F_1(a,b,a+b+1,x),$$
$$g(x)=-\Big(2(ab/a+b) _2F_1(a,b,a+b+1,x)+\,_2F_1(a,b,a+b,x)\Big),$$
and $$\Delta(x)=(g(x))^2-8h(x)\,_2F_1(a,b,a+b,x).$$
\begin{theorem}\label{t1} For $a,b\geq 0$ and $c\in\Bbb R$, we define
$$f_{a,b,c}(x)=\frac{\,_2F_1(a,b,a+b,x)}{c-\log (1-x)}.$$
Assume $a+b\leq 1$, then the function $f_{a,b,c}$ is convex on $[0,1)$ if and only if $c\geq \alpha_0$ and is concave on $[0,1)$ if and only if
$c\in[\delta_-,\delta_+]$, where
$\alpha_0=\max_{x\in[0,1)}\varphi_+(x)$, $\delta_-=\max_{x\in[0,1)}\varphi_-(x)$ and $\delta_+=\min_{x\in[0,1)}\varphi_+(x)$.
\end{theorem}
Observe that $\delta_+>\delta_-$. We will see that equality holds if and only if $a=b=1/2$, in that case $\delta_+=8/3$.


\section{Lemmas}
In this section, we collect some results which are needed to prove Theorem \ref{t2} and Theorem \ref{t1}. The first lemma offers three basic properties of the hypergeometric function ${}_2F_1$ (see [\cite{olv}, 15.4.20, 15.4.21, 15.5.1,
15.8.1, 15.8.11], Abramowitz 15.3.10) and \cite{ask}.

\begin{lemma}\label{lem1} For $x\in(-1,1)$
$$\frac{d}{dx}{}_2F_1(a,b;c,x)=\frac{ab}c{}_2F_1(a+1,b+1;c+1,x),$$
$${}_2F_1(a,b;c,x)=(1-x)^{c-b-a}{}_2F_1(c-a,c-b;c,x),$$
$${}_2F_1(a,b;a+b,x)=\frac{1}{B(a,b)}\sum_{n=0}^\infty\frac{(a)_n(b)_n}{(n!)^2}\Big(2\psi(n+1)-\psi(a+n)-\psi(b+n)-\log(1-x)\Big)(1-x)^n,$$
for $x\in(0,2)$, and if $c>a+b$,
$${}_2F_1(a,b;c,1)=\frac{\Gamma(c)\Gamma(c-b-a)}{\Gamma(c-b)\Gamma(c-a)}.$$
\end{lemma}
\begin{lemma}\label{lem2} Let $A(x)=\sum_{k=0}^\infty a_kx^k$ and $B(x)=\sum_{k=0}^\infty a_kx^k$ be two reals power series converging on $(-r,r)$, $(r>0)$ with $b_k>0$ for all $k$, if the non-constant sequence $(a_k/b_k)_k$ is increasing (resp. decreasing) for all $k$, then the function $x\mapsto A(x)/B(x)$ is strictly increasing (resp. decreasing) on $(0,r)$.
\end{lemma}
\begin{lemma}\label{lem3} Let $A(x)=\sum_{k=0}^\infty a_kx^k$ and $B(x)=\sum_{k=0}^\infty a_kx^k$ be two reals power series converging on $(-r,r)$, $(r>0)$ with $b_k>0$ for all $k$. Suppose that for certain $m\in\Bbb N$ the non-constant sequence $(a_k/b_k)_k$ is increasing (resp. decreasing) for all $0\leq k\leq m$, and
decreasing (resp. increasing) for $k\geq m$. Then the function $x\mapsto A(x)/B(x)$ is strictly increasing (resp. decreasing) on $(0,r)$ if and only if
$H_{A,B}(r^-)\geq 0$ (resp. $H_{A,B}(r^-)\leq 0$), where
$H_{A,B}(x)=(A'(x)/B'(x))B(x)-A(x)$.
\end{lemma}
See for instance \cite{bier}, (see also [\cite{ponn}, Lemma 2.1], \cite{hei}).
\begin{lemma}\label{lem5} Let $\alpha,\beta\in\Bbb R$ and $f,g$ be two continuous functions on $[\alpha,\beta]$ and differentiable on $(\alpha,\beta)$ and $g'(x)\neq 0$ for all $x\in(\alpha,\beta)$. If the ration $f'(x)/g'(x)$ is increasing (respectively decreasing) on $(\alpha,\beta)$ then so are the ratios
$$\frac{f(x)-f(a)}{g(x)-g(a)},\quad \frac{f(x)-f(b)}{g(x)-g(b)}.$$
If  $f'(x)/g'(x)$ is strictly monotone, then the monotonicity in the conclusion is also strict.
\end{lemma}
\begin{lemma}\label{l1} For $a,b\geq 0$ and $a+b\leq 1$, we have
\begin{enumerate}
\item $4ab\leq a+b$.
 \item The function $x\mapsto S(a,b,x)=1-\frac{a^2b^2}{(a+b+1)(a+b)}-\frac{ab(a+1)(b+1)(a+b+ab+2)}{(a+b)(a+b+1)(a+b+2)}\,x$ is negative for all $x\geq c(a,b)$
 where, $$c(a,b)=\frac{(a+b)(a+b-2 a b)(a+b+1)}{a b ((a+b+1)(a+b-2 a b)+ab(a+b))}.$$
     \end{enumerate}
\end{lemma}
{\bf Proof.} 1) The case $a=0$ or $b=0$ is trivial. For the case $a\neq 0$ and $b\neq 0$, the inequality $4ab\leq a+b$ is equivalent to
$1/a+1/b\geq 4$. Therefore, if $1/a+1/1-a \geq 4$ then $1/a+1/b\geq 4$. Let
$h(a)=1/a+1/1-a$, then $h'(a)=1/(1-a)^2 -1/a^2$. Then $h$ decreases on $(0,1/2]$ and increases on $[1/2,1)$ and $h(1/2)=4$.

2) 
We have, $$S(a,b,x)\leq S(a,b,c(a,b))=-ab\frac{P(a,b)}{Q(a,b)}.$$
Where,
$$ Q(a,b)=(a+b+2)((a+b+1)(a+b-2 a b)+ab(a+b))(a+b+1)(a+b)$$
$$P(a,b)=Q_0(b)+Q_1(b)a+Q_2(b)a^2+Q_3(b)a^3+Q_4(b)a^4,$$
$$Q_0(b)=b^2(1+b)^2,\quad Q_1(b)=2(1 - b) b (1 + b)^2,$$
$$Q_2(b)=1 + 2 b - 10 b^2 - 12 b^3 - 3 b^4,\quad Q_3(b)=2(1 -b - 6 b^2 -3 b^3),$$
and
$$Q_4(b)=(1+b)(1-3b).$$
By the first item we have $Q(a,b)>0$ for all $a+b\leq 1$ and $a,b>0$ and $Q_0(b)>0$, $Q_1(b)>0$.

From the expression of $S(a,b,c(a,b))$ one see that it is a symmetric function on the variables $a,b$ and by the condition $a+b\leq 1$, at least one of the variables $a$ or $b$ is in $[0,1/2]$. So, one may assume $b\in[0,1/2]$ and $a\in[0,1-b]$. Differentiate yields
$$\frac{\partial^4}{\partial^4 a}P(a,b)=24 Q_4(b),$$
which is positive for $b\in(0,1/3]$ and negative for $b\in[1/3,1/2]$.

{\bf Case 1:} $b\in[1/3,1/2].$
Firstly, we have $Q'_3(b)=-2 (1 + 12 b + 9 b^2)<0$ and $Q_3(b)<Q_3(1/3)=-2/9$
then, for $a\in(0,1-b)$, $$0>6Q_3(b)>\frac{\partial^3}{\partial^3 a}P(a,b).$$
Therefore, $a\mapsto \frac{\partial^2}{\partial^2 a}P(a,b)$ decreases on $[0,1-b]$ for all $b\in[1/3,1/2].$
Moreover, $$\frac{\partial^2}{\partial^2 a}P(a,b)=2Q_2(b)+6Q_3(b)a+12Q_4(b)a^2,$$
and $\frac{\partial^2}{\partial^2 a}P(0,b)=2Q_2(b)>0$ and $\frac{\partial^2}{\partial^2 a}P(1-b,b)<0$ for $b\in[1/3,1/2]$.
So, there is $a_0\in(0,1-b)$, such that the  function $a\mapsto \frac{\partial}{\partial a}P(a,b)$ increases on $[0,a_0]$ and decreases on $[a_0,1-b]$.

 On the other hand,
$$\frac{\partial}{\partial a}P(a,b)=Q_1(b)+2Q_2(b)a+3Q_3(b)a^2+4Q_4(b)a^3,$$
and $$\frac{\partial}{\partial a}P(0,b)=Q_1(b)>0, \frac{\partial}{\partial a}P(1-b,b)=Q_1(b)+2Q_2(b)(1-b)+3Q_3(b)(1-b)^2+4Q_4(b)(1-b)^3,$$
then $$P(a,b)\geq\min(P(0,b),P(1-b,b),$$
Besides, $P(0,b)=Q_0(b)>0,$ and $$P(1-b,b)=Q_0(b)+(1-b)Q_1(b)+(1-b)^2Q_2(b)+(1-b)^3Q_3(b)+(1-b)^4Q_4(b),$$
hence, $$\theta(b):=P(1-b,b)=-13 b^4+26 b^3-b^2-12 b+4.$$
Since, $\theta''(b)=-2 (1 - 78 b + 78 b^2)$ and $\theta''(b)=0$ on $[0,1/2]$ for $b_0=1/78 (39 - \sqrt{1443})$, and is negative $[0,b_0]$ and positive on $[b_0,1/2]$. Moreover, $\theta'(b)=-2 (6 + b - 39 b^2 + 26 b^3)$ with $\theta'(0)=-12$, $\theta'(1/2)=0$. Thus, $\theta(b)$ decreases on $[0,1/2]$ and we have
$\theta(1/2)=3/16$. Therefore, $P(1-b,b)>0$ for all $b\in[0,1/2]$ and then $P(a,b)>0$ for all $b\in[1/3,1/2]$ and $a\in[0,1-b]$.

{\bf Case 2:} Assume $b\in[0,1/3]$.

It is easy to see that $Q_0$, $Q_1$ and $Q_4$ are positive on $[0,1/3)$. Moreover,
$Q''_2(b)=-4 (5 + 18 b + 9 b^2)<0$ and $Q'_2(b)=-2 (-1 + 10 b + 18 b^2 + 6 b^3)$ with $Q'_2(0)=2$ and $Q'_2(1/3)=-82/9$. Therefore, $Q_2(b)\geq\min(Q_2(0),Q_2(1/3))=2/27$, whence, $Q_2(b)>0$.

By straightforward computation, one see that there is a unique $b_0\in(0,1/3)$ such that $Q_3(b)>0$ on $(0,b_0)$ and $Q_3(b)<0$ on $(b_0,1/3)$. In that case $P(a,b)>0$ for all $b\in(0,b_0)$.

Assume $b\in[b_0,1/3]$:

Since, $Q_4(b)>0$ on $[0,1/3)$, then $a\mapsto \frac{\partial^3}{\partial^3 a}P(a,b)$ increases on $[0,1-b]$ with $\frac{\partial^3}{\partial^3 a}P(0,b)=6Q_3(b)<0,$ on $[b_0,1/3]$ and $\frac{\partial^3}{\partial^3 a}P(1-b,b)>0$. Whence, there is a unique $a_0\in(0,1-b)$ such that $a\mapsto \frac{\partial^2}{\partial^2 a}P(a,b)$ decreases on $(0,a_0]$ and increases on $(a_0,1-b]$. Furthermore,
$$\frac{\partial^2}{\partial^2 a}P(a,b)=2Q_2(b)+6Q_3(b)a+12Q_4(b)a^2.$$
If $\frac{\partial^2}{\partial^2 a}P(a_0,b)>0$, then $\frac{\partial}{\partial a}P(a,b)$ increases and since $\frac{\partial}{\partial a}P(0,b)=Q_1(b)>0$ and then $P(a,b)>P(0,b)=Q_0(b)>0$.

If $\frac{\partial^2}{\partial^2 a}P(a_0,b)<0$, then, there are $a_1\in(0,a_0)$ and $a_2\in(a_0,1-b)$ such that $a\mapsto \frac{\partial}{\partial a}P(a,b)$  increases on $(0,a_1)$ and on $(a_2,1-b)$ and decreases on $(a_1,a_2)$. Moreover,
 $$\frac{\partial^2}{\partial^2 a}P(a_2,b)=2Q_2(b)+6Q_3(b)a_2+12Q_4(b)a_2^2=0.$$
 So,$$\frac{\partial}{\partial a}P(a_2,b)= Q_1(b)+Q_2(b)a_2-2Q_4(b)a_2^3\geq Q_1(b)-2Q_4(b)(1-b)^3:=\theta_1(b).$$
 We have $\theta_1(b)=2(1 - 4 b + 8 b^2 - 3 b^3)(1-b^2)$ and $(\theta_1(b)/2(1-b^2))'=-4+16b-9b^2$, which is negative on $(0, 2/9 (4 - \sqrt 7))$ and positive on $(2/9 (4 - \sqrt 7) ,1/3)$. Then $1 - 4 b + 8 b^2 - 3 b^3\geq (1/243) (403 - 112\sqrt 7)>0,$  and $\theta_1(b)>0$ for $b\in(b_0,1/3)$. Whence, $\frac{\partial}{\partial a}P(a_2,b)>0$. So, $$\frac{\partial}{\partial a}P(a,b)\geq\min(\frac{\partial}{\partial a}P(a_2,b),\frac{\partial}{\partial a}P(0,b))>0.$$
 Therefore, the function
 $a\mapsto P(a,b)$ is strictly increasing on $(0,1-b)$ and $P(a,b)>P(0,b)=Q_0(b)>0$.

This completes the proof.

\begin{lemma}\label{i} Let $a,b\geq 0$ and $a+b\leq 1$.
Then, the function $h(x)$ is strictly increasing from $[0,1)$ onto $\Big[\frac{a(a+1)b(b+1)}{(a+b)(a+b+1)},\frac1{B(a,b)}\Big]$ and the function $\Delta(x)$ is strictly increasing from $(0,1)$ onto $(0,+\infty)$.

 The function $\varphi_+(x)$ initially defined on $[0,1)$, is extended to a continuous function on $[0,1]$ with $\varphi_+(1)=R(a,b)$. It admits a maximum $\alpha_0$ on $[0,1)$.

 The function $\varphi_-(x)$ is continuous on $[0,1)$ with $\displaystyle\lim_{x\to 1^-}w_-(x)=-\infty$.
\end{lemma}
{\bf Proof.} We have
$$\begin {aligned}\Delta(x)&=4(\frac{ab}{a+b})^2( _2F_1(a,b,a+b+1,x))^2+\,_2F_1(a,b,a+b,x)^2\\
&-4\frac{ab}{a+b}\,_2F_1(a,b,a+b,x) _2F_1(a,b,a+b+1,x)\\
&-8\frac{(ab)^2}{(a+b)(a+b+1)}\,_2F_1(a+1,b+1,a+b+2,x)F_1(a,b,a+b,x)(1-x)\end{aligned}$$
Differentiate yields
$$\begin{aligned}\Delta'(x)&=8(\frac{ab}{a+b})^2\frac{ab}{a+b+1}\,_2F_1(a,b,a+b+1,x)\,_2F_1(a+1,b+1,a+b+2,x)\\
&+2\frac{ab}{a+b}\,_2F_1(a+1,b+1,a+b+1,x)\,_2F_1(a,b,a+b,x)\\
&-4(\frac{ab}{a+b})^2\,_2F_1(a+1,b+1,a+b+1,x) _2F_1(a,b,a+b+1,x)\\
&-4\frac{(ab)^2}{(a+b)(a+b+1)}\,_2F_1(a,b,a+b,x) _2F_1(a+1,b+1,a+b+2,x)\\
&-8\frac{(ab)^2(a+1)(b+1)}{(a+b)(a+b+1)(a+b+2)}\,_2F_1(a+1,b+1,a+b+3,x)F_1(a,b,a+b,x\\
&-8\frac{(ab)^3}{(a+b)^2(a+b+1)}\,_2F_1(a+1,b+1,a+b+2,x)F_1(a,b,a+b+1,x)\\
&+8\frac{(ab)^2}{(a+b)(a+b+1)}\,_2F_1(a+1,b+1,a+b+2,x)F_1(a,b,a+b,x)\end{aligned}$$
$$\begin{aligned}\frac{a+b}{2ab}\Delta'(x)&=\,_2F_1(a+1,b+1,a+b+1,x)\,_2F_1(a,b,a+b,x)\\
&-2\frac{ab}{a+b}\,_2F_1(a+1,b+1,a+b+1,x) _2F_1(a,b,a+b+1,x)\\
&-4\frac{(ab)(a+1)(b+1)}{(a+b+1)(a+b+2)}\,_2F_1(a+1,b+1,a+b+3,x)F_1(a,b,a+b,x)\\
&+2\frac{ab}{(a+b+1)}\,_2F_1(a+1,b+1,a+b+2,x)\,_2F_1(a,b,a+b,x)\end{aligned}$$
Set $$\Delta_1(x)=\frac12\,_2F_1(a,b,a+b,x)-\frac{2ab}{a+b}\,_2F_1(a,b,a+b+1,x),$$ and
$$\Delta_2(x)=\frac12\,_2F_1(a+1,b+1,a+b+1,x)+\frac{2ab}{(a+b+1)}\,_2F_1(a+1,b+1,a+b+2,x)$$
$$-4\frac{(ab)(a+1)(b+1)}{(a+b+1)(a+b+2)}\,_2F_1(a+1,b+1,a+b+3,x).$$
It follows that
$$\frac{a+b}{2ab}\Delta'(x)=\Delta_1(x)\,_2F_1(a+1,b+1,a+b+1,x)+\Delta_2(x)\,_2F_1(a,b,a+b,x).$$
Using the formula $(a+1)_n=((a+n)/a)(a)_n$, then
$$\Delta_1(x)=\frac12\sum_{n=0}^\infty\frac{a+b-4ab+n}{a+b+n}\frac{(a)_n(b)_n}{n!(a+b)_n}x^n.$$
Using Lemma \ref{l1}, we get $\Delta_1(x)\geq 0$.
On the other hand
$$\Delta_2(x)=\frac12\sum_{n=0}^\infty\alpha_n\frac{(a+1)_n(b+1)_n}{n!(a+b+1)_n}x^n,$$
where $$\alpha_n=1+\frac{4ab}{a+b+1+n}-\frac{8ab(a+1)(b+1)}{(a+b+1+n)(a+b+2+n)}.$$
From Lemma \ref{l1}, we have $ab\leq 1/4$, so
$$\begin{aligned}1-\frac{2(a+1)(b+1)}{a+b+2+n}&=\frac{a+b+2+n-2(a+1)(b+1)}{a+b+2+n}\\&=-\frac{a+b+2ab-n}{a+b+2+n}\\&\geq -\frac{a+b+1/2-n}{a+b+2+n}\geq -1,\end{aligned}$$
 then, $$\begin{aligned}\alpha_n&\geq 1+\frac{4ab}{a+b+1}(1-\frac{2(a+1)(b+1)}{a+b+2})\\
&\geq \frac{a+b-4ab+1}{a+b+1}>0.\end{aligned}$$
and $\Delta_2(x)> 0$. Hence, $\Delta'(x)>0$ and $\Delta(x)$ is strictly increasing. Moreover,
$$\Delta(0)=1-4\frac{ab}{a+b}+4(\frac{ab}{a+b})^2\frac{1-a-b}{a+b+1}> 0,$$
equality holds if and only if $a=b=1/2$.

On the first hand, we have $\displaystyle\lim_{x\to 1^-}\,_2F_1(a,b,a+b,x)=+\infty$ and for $x$ close to $1$, $g(x)\sim \,_2F_1(a,b,a+b,x)$.
Hence,
$$\Delta(x)=-g(x)(1- \frac{4h(x)\,_2F_1(a,b,a+b,x)}{g(x)^2}+\frac18(8\frac{h(x)\,_2F_1(a,b,a+b,x)}{g(x)^2})^2+o(\frac{\,_2F_1(a,b,a+b,x)}{(g(x))^2})),$$
and $$\omega_+(x)=-\frac{g(x)}{h(x)}+\frac{2\,_2F_1(a,b,a+b,x)}{g(x)}-\frac{g(x)}{16h(x)}\Big(\frac{h(x)\,_2F_1(a,b,a+b,x)}{g(x)^2}\Big)^2
+o((\frac{\,_2F_1(a,b,a+b,x)}{g(x)})^2).$$
Since, $h(1)=1/B(a,b)$, then
$$\begin{aligned}\varphi_+(x)&=\log(1-x)-B(a,b)g(x)+\frac{2\,_2F_1(a,b,a+b,x)}{g(x)}-\frac{g(x)}{16h(x)}
\Big(\frac{h(x)\,_2F_1(a,b,a+b,x)}{g(x)^2}\Big)^2\\
&+o((\frac{\,_2F_1(a,b,a+b,x)}{g(x)})^2).\end{aligned}$$
Therefore, $$\varphi_+(x)=\log(1-x)-B(a,b)g(x)+2
+o((\frac1{\,_2F_1(a,b,a+b,x)})),$$
$$\varphi_+(x)=\log(1-x)-B(a,b)g(x)-2+o(\frac1{\,_2F_1(a,b,a+b,x)}).$$
As, $$g(x)=-(2(ab/a+b) _2F_1(a,b,a+b+1,x)+\,_2F_1(a,b,a+b,x)),$$
we get $$\varphi_+(x)=\log(1-x)+B(a,b)\,_2F_1(a,b,a+b,x)+o(\frac1{\,_2F_1(a,b,a+b,x)}),$$
and from Lemma \ref{lem1}, we have
$$\lim_{x\to 1^-}\varphi_+(x)=R(a,b),$$
and $$\varphi_+(0)=\frac{-g(0)+\sqrt{\Delta(0)}}{2h(0)}=(a+b+1)\frac{a+b+2ab+(a+b)\sqrt{\Delta(0)}}{2ab(a+b+ab+1)}>0.$$
So, the function $\varphi_+$ is bounded on $[0,1]$. Moreover,
$\lim_{x\to 1^-}\varphi_-(x)=-\infty$ and
$$\varphi_-(0)=(a+b+1)\frac{a+b+2ab-(a+b)\sqrt{\Delta(0)}}{2ab(a+b+ab+1)}.$$

For $a=b=1/2$, we get $R(1/2,1/2)=\log 16$.
Which completes the proof of the lemma.
\begin{proposition}\label{pro4} For every $b>0$, define the functions
$$H_b(x)=B(x,b)-R(x,b),$$
$$G_b(x)=B(x,b)-\frac1x.$$
\begin{enumerate}
\item The function $H'_b$ is strictly completely monotonic on $(0,+\infty)$ if and only if $b>0$.

\item The function $G_b$ is strictly completely monotonic on $(0,+\infty)$ if and only if $b\in(0,1)$.

 \end{enumerate}
 As consequence, for all $a>0$ and $b\in(0,1)$
 $$R(a,b)<B(a,b)<\min(U(a,b),U(b,a))<\frac1a+\frac1b,$$
 where $U(a,b)=-\gamma+1/a-\psi(b)$.
\end{proposition}
{\bf Proof.}
1) Since, $$H_b(x)=\gamma+\psi(b)+\int_0^1t^{x-1}(1-t)^{b-1}dt+\int_0^1\frac{1-t^{x-1}}{1-t}dt.$$
Hence,
$$H'_b(x)=\int_0^1t^{x-1}(1-t)^{b-1}\log (t)dt-\int_0^1\frac{t^{x-1}\log(t)}{1-t}dt.$$
Hence, $$H'_b(x)=\int_0^1t^{x-1}((1-t)^{b}-1)\frac{\log (t)}{1-t}dt.$$
Therefore, for all $n\geq 0$,
 $$(-1)^nH^{(n+1)}_b(x)=\int_0^1t^{x-1}(1-(1-t)^{b})\frac{(-\log (t))^n}{1-t}dt,$$
 and $H_b'$ is strictly completely monotonic on $(0,+\infty)$ if and only if $b>0$.
 So, the function $H_b$ is strictly increasing on $(0,+\infty)$. Further, for $x$ close to $0$,
$\log\Gamma(x)=-\log x-\gamma x+o(x),$ $\psi(x)=-\frac1x+\frac{\pi^2}2 x+o(x),$ and $\log\Gamma(x+b)=\log\Gamma(b)+\psi(b)x+\frac12\psi'(b)x^2+o(x)$.
Then, $$H_b(x)=\frac{1}{4}\left(2 \psi(b)^2+4 \gamma\psi(b)-2 \psi'(b)+2 \gamma ^2+\pi ^2\right)x+o(x),$$
and $\displaystyle\lim_{x\to 0^+}H_b(x)=0$. This proves the first part of the inequality.

2) We have $$G_b(x)=\int_0^1 t^{x-1}((1-t)^{b-1}-1) dt.$$
So, $G_b$ is strictly completely monotonic if and only if $b-1<0$.

Since, $\lim_{x\to +\infty} G_b(x)=0$ and $\lim_{x\to 0}G_b(x)=-\gamma-\psi(b)$, which gives the second part of the inequality. Moreover, the function $b\mapsto-\gamma-\psi(1+b)$ is strictly decreasing and negative, then $U(a,b)=-\gamma+1/a-\psi(b)=1/a+1/b-\gamma-\psi(1+b)<1/a+1/b$.
\begin{lemma}\label{l0} For $a,b> 0$, the function $\varphi_{a,b,c}$ is strictly increasing from $(0,1)$ onto $(c, B(a,b))$ if and only if $c\leq R(a,b)$ and is strictly decreasing on $(0,1)$ if and only if $c\geq 1/a+1/b$.
\end{lemma}
For $a=b=1/2$ we retrieve, the results of [Anderson et al], in that case $c\leq4\log 2$ for increasing and $c\geq 4$ for the decreasing.

{\bf Proof.} Differentiate yields
$$\varphi'_{a,b}(x)=\frac{\,_2F_1(a,b,a+b,x)-\frac{ab}{a+b}\,_2F_1(a,b,a+b+1,x)(c-\log(1-x))}{(1-x)\,_2F_1(a,b,a+b,x)^2}=\frac{u(x)}{v(x)},$$
and
$$u'(x)=-\frac{(ab)^2}{(a+b)(a+b+1)}\,_2F_1(a+1,b+1,a+b+2,x)(c-\log(1-x)).$$
Therefore, $u(x)$ is strictly decreasing. Moreover,
\begin{equation}\label{lim}\lim_{x\to 1^-} u(x)=\frac{R(a,b)-c}{B(a,b)},\end{equation}
and $u(0)=(ab/(a+b))(1/a+1/b-c)$. Which proves the lemma.
\begin{lemma} For $a,b\geq 0$, let
$$A(x)=(1-x)\,_2{F}_1(a,b;a+b;x)\,_2{F}_1(a,b;a+b+1;x),$$
then $A'(x)<0$ for all $x\in(0,1)$.
\end{lemma}
{\bf Proof.} Differentiate yields
$$\begin{aligned}A'(x)&=\frac{ab}{a+b}\,_2{F}_1(a,b;a+b+1;x)^2+
\frac{ab(1-x)}{a+b+1}\,_2{F}_1(a,b;a+b;x)\,_2{F}_1(a+1,b+1;a+b+2;x)\\
&-\,_2{F}_1(a,b;a+b;x)\,_2{F}_1(a,b;a+b+1;x).\end{aligned}$$
$$\,_2{F}_1(a,b;a+b+1;x)^2\leq(\frac{\Gamma(a+b+1)}{\Gamma(a+1)\Gamma(b+1)})^2,$$
$$\,_2{F}_1(a,b;a+b;x)\,_2{F}_1(a,b;a+b+1;x)\geq 1,$$ then
$$A'(x)\leq -1+\frac{ab}{a+b}(\frac{\Gamma(a+b+1)}{\Gamma(a+1)\Gamma(b+1)})^2+\frac{ab}{a+b+1}L(x),$$
where $$L(x)=(1-x)\,_2{F}_1(a,b;a+b;x)\,_2{F}_1(a+1,b+1;a+b+2;x).$$
Moreover,
$$\begin{aligned}L'(x)&=\frac{(a+1)(b+1)}{a+b+2}\,_2{F}_1(a,b;a+b;x)\,_2{F}_1(a+1,b+1;a+b+3;x)\\
&+\frac{ab}{a+b}\,_2{F}_1(a,b;a+b+1;x)\,_2{F}_1(a+1,b+1;a+b+2;x)\\&-\,_2{F}_1(a,b;a+b;x)\,_2{F}_1(a+1,b+1;a+b+2;x),
\end{aligned}$$
$$\begin{aligned}L'(x)&=M\,_2{F}_1(a,b;a+b;x)\Big(\,_2{F}_1(a+1,b+1;a+b+3;x)-\,_2{F}_1(a+1,b+1;a+b+2;x)\Big)\\
&+\frac{ab}{a+b}\,_2{F}_1(a+1,b+1;a+b+2;x)\Big(\,_2{F}_1(a,b;a+b+1;x)-\,_2{F}_1(a,b;a+b;x)\Big)\\
&-M'\,_2{F}_1(a,b;a+b;x)\,_2{F}_1(a+1,b+1;a+b+2;x),\end{aligned}$$
where $M=\frac{(a+1)(b+1)}{a+b+2}$ and $M'=\frac{(1-ab)(a+b)+ab(a+b+2)}{(a+b+2)(a+b)}$.

Since, $\,_2{F}_1(a,b;c+1;x)<\,_2{F}_1(a,b;c;x)$ for ($a,b,c>0$) and $a\leq 1$, $b\leq 1$, it follows that
$L'(x)\leq 0$ and hence $L(x)<L(0)=1$. This implies that
$$A'(x)\leq -1+\frac{ab}{a+b}(\frac{\Gamma(a+b+1)}{\Gamma(a+1)\Gamma(b+1)})^2+\frac{ab}{a+b+1}:=\theta(a,b),$$
let $m(a)=\log\Gamma(a+b+1)-\log\Gamma(a+1)-\log\Gamma(b+1),$ then
$m'(a)=\psi(a+b+1)-\psi(a+1)$, moreover the function $\psi$ increases then $m'(a)\geq 0$ and the function $m(a)$ increases. Then the function $a\mapsto \Gamma(a+b+1)/\Gamma(a+1)\Gamma(b+1)$ increases and is positive. Furthermore, the functions $a\mapsto a b/(a+b)$ and $a\mapsto ab/a+b+1$ increases. Besides, the function $\theta(a,b)$ is symmetric, then
 $\theta(a,b)\leq\theta(\frac12,1)=-1/20.$
 \begin{lemma} Let $a,b\geq 0$, $a+b\leq 1$ and
 $$B_{a,b}(x)=\,_2{F}_1(a,b;a+b;x)-\frac{a(a+1)b(b+1)}{(a+b)(a+b+1)}\,_2{F}_1(a,b;a+b+2;x)
 (c-\log(1-x)).$$
 Then, $(-B'_{a,b}(x))^{(n)}>0$ for all $n\geq 0$ and all $x\in(0,1)$.
 \end{lemma}
 {\bf Proof.} We have
 $$-B'_{a,b}(x)=\frac{ab}{a+b}\frac1{1-x}(\frac{(a+1)(b+1)}{(a+b+1)}\,_2{F}_1(a,b;a+b+2;x)-\,_2{F}_1(a,b;a+b+1;x))$$
 $$+\frac{a^2(a+1)b^2(b+1)}{(a+b)(a+b+1)(a+b+2)}\,_2{F}_1(a,b;a+b+3;x)
 (c-\log(1-x)$$
 Let us set, $$\frac{(a+1)(b+1)}{(a+b+1)}\,_2{F}_1(a,b;a+b+2;x)-\,_2{F}_1(a,b;a+b+1;x):=C(x),$$
then $$C(x)=\sum_{n=0}^\infty c_nx^n,$$
$c_0=a^2b^2/(a+b)(a+b+1)$ and $c_n=\frac{ab(a)_n(b)_n}{(a+b)(a+b+1)_n n!}(\frac{(a+1)(b+1)}{(a+b+1+n)}-1)$ for $n\geq 1$.
 Since, $\frac{(a+1)(b+1)}{(a+b+1+n)}-1=-(n-ab)/(a+b+1+n)$, then $c_n<0$ for $n\geq 1$. Moreover,
 $$\frac{C(x)}{1-x}=\sum_{n=0}^\infty \alpha_n x^n,$$
 and
 $$\begin{aligned}\alpha_n&=\sum_{k=0}^n c_k\geq \sum_{k=0}^\infty c_k=C(1)\\
 &=\frac{(a+1)(b+1)}{(a+b+1)}\,_2{F}_1(a,b;a+b+2;1)-\,_2{F}_1(a,b;a+b+1;1)\\
 &=\frac{(a+1)(b+1)\Gamma(a+b+2)}{(a+b+1)\Gamma(b+2)\Gamma(a+2)}-\frac{\Gamma(a+b+1)}{\Gamma(a+1)\Gamma(b+1)}=0.\end{aligned}$$
 Since,
  $$-B'_{a,b}(x)=\frac{C(x)}{1-x}+K(x),$$
  where $$K(x)=\frac{a^2(a+1)b^2(b+1)}{(a+b)(a+b+1)(a+b+2)}\,_2{F}_1(a,b;a+b+3;x)
 (c-\log(1-x)).$$
 Moreover, the function $K(x)$ is absolutely monotonic, then $-B'_{a,b}(x)$ is absolutely monotonic.
 \begin{lemma} For $a,b\geq 0$ and $a+b\leq 1$, we set
 $$P_{a,b}(x)=\frac{B'_{a,b}(x)A(x)}{A'(x)}-B_{a,b}(x),$$
 then $P_{a,b}(x)>0$ for all $x\in(0,1)$.
 \end{lemma}
 {\bf Proof.}
 We have
 $$A'(x)>-\,_2{F}_1(a,b;a+b;x)\,_2{F}_1(a,b;a+b+1;x)=-\frac{A(x)}{1-x},$$
 therefore,
 $$-\frac{A(x)}{A'(x)}\geq (1-x),$$
 and
  $$\frac{B'_{a,b}(x) A(x)}{A'(x)}\geq  -(1-x)B'_{a,b}(x).$$
  We set $S_{a,b}(x)=B_{a,b}(x) +(1-x)B'_{a,b}(x)$. We have $S'_{a,b}(x)=(1-x)B''_{a,b}(x)<0$, therefore, $S_{a,b}(x)<B_{a,b}(0)+B'_{a,b}(0)$
Since, $$\begin{aligned}B_{a,b}(0)+B'_{a,b}(0)&=1-\frac{a^2b^2}{(a+b)(a+b+1)}-c\Big(\frac{ab(a+1)(b+1)(a+b+2)+a^2b^2(a+1)(b+1)}{(a+b)(a+b+1)(a+b+2)}
\Big)\\&=S(a,b,c),\end{aligned}$$
 where $S(a,b,c)$ is the function defined in Lemma \ref{l1}. It has been shown that $S(a,b,c)<0$ for $c\geq c(a,b)$
  then, $P_{a,b}(x)\geq -S_{a,b}(x)>0$.

\section{Proofs of theorems}

\subsection{Proof of Theorem \ref{t2}}
1) Assume that $a+b\leq 1$ and $c\leq R(a,b)$. From Lemma \ref{lem1}, we have
$$\varphi'_{a,b,c}(x)=\frac{\,_2F_1(a,b,a+b,x)-\frac{ab}{a+b}\,_2F_1(a,b,a+b+1,x)(c-\log(1-x))}{(1-x)\,_2F_1(a,b,a+b,x)^2}=\frac{u(x)}{v(x)}.$$
$$u'(x)=-\frac{(ab)^2}{(a+b)(a+b+1)}\,_2F_1(a+1,b+1,a+b+2,x)(c-\log(1-x)),$$
\begin{equation}\label{mon}v'(x)=\frac{2ab}{a+b}\,_2F_1(a,b,a+b,x)\,_2F_1(a,b,a+b+1,x)-\,_2F_1(a,b,a+b,x)^2.\end{equation}
Then
$$\frac{u'(x)}{v'(x)}=\frac{(ab)^2}{(a+b)(a+b+1)}\frac{\,_2F_1(a+1,b+1,a+b+2,x)}{\,_2F_1(a,b,a+b,x)-\frac{2ab}{a+b}
\,_2F_1(a,b,a+b+1,x)}\varphi_{a,b,c}(x).$$
Let $$\Theta(x)=\frac{\,_2F_1(a+1,b+1,a+b+2,x)}{\,_2F_1(a,b,a+b,x)-\frac{2ab}{a+b}
\,_2F_1(a,b,a+b+1,x)}.$$
Direct computation gives
$$\Theta(x)=\frac{\sum_{n=0}^\infty a_nx^n}{\sum_{n=0}^\infty b_nx^n},$$
where, $$a_n=(a+1)_n(b+1)_n/(n!(a+b+2)_n),$$ and
$$b_n=(a)_n(b)_n/(n!(a+b)_n))-(2ab/(a+b))((a)_n(b)_n/(n!(a+b+1)_n)).$$
 Using the formula $(x+1)_n=((x+n)/x)(x)_n$, we get
$$\frac{a_n}{b_n}=\frac{(a+b+1)(a+b)}{ab}\frac{(a+n)(b+n)}{(a+b+1+n)(a+b+n-2ab)}:=u_n(a,b).$$
Let us write
$$u_n(a,b)=\frac{(a+b+1)(a+b)}{ab}\frac{1}{(1+\frac{b+1}{a+n})(1+\frac{a(1-2b)}{b+n})}.$$
If $b\in(0,1/2]$ the sequence $u_n(a,b)$ is strictly increasing. If $a\in(0,1/2]$ then by using
the symmetry $u_n(a,b)=u_n(b,a)$ we deduce that the sequence $u_n(a,b)$ is strictly increasing and so is $a_n/b_n$.

Therefore, the function $u'(x)/v'(x)$ is strictly increasing and from Lemma \ref{l0}, one deduces that $\displaystyle\frac{u'(x)}{v'(x)}=\frac{(ab)^2}{(a+b)(a+b+1)}\Theta(x)\varphi_{a,b,c}(x)$ is strictly increasing too. From Lemma \ref{lem5}, we get that the function
$\displaystyle\frac{u(x)-u(1^-)}{v(x)-v(1^-)}$ is strictly increasing. Besides, from equation \eqref{lim}, we have $u(1^-)=\frac{-c+R(a,b)}{B(a,b)}\geq 0$ and $v(1^-)=0$ and by equation \eqref{mon} and the fact that $2ab\leq a+b$, the function $v(x)$ is strictly decreasing and positive. Therefore, the function $\displaystyle\frac{u(x)}{v(x)}=\frac{u(x)-u(1^-)}{v(x)}+\frac{u(1^-)}{v(x)}$ is strictly increasing.
 It then follows that, the function $\varphi'_{a,b,c}(x)$ is strictly increasing and the function $\varphi_{a,b,c}(x)$ is strictly convex on $(0,1)$.

2) Next, assume that $\varphi_{a,b,c}$ is strictly convex, then $\varphi'_{a,b,c}$ is strictly increasing. Recall that
$$\varphi'_{a,b}(x)=\frac{\,_2F_1(a,b,a+b,x)-\frac{ab}{a+b}\,_2F_1(a,b,a+b+1,x)(c-\log(1-x))}
{(1-x)\,_2F_1(a,b,a+b,x)^2}.$$
From Lemma \ref{lem1}, we have $\,_2F_1(a,b,a+b+1,1)=\Gamma(a+b+1)/\Gamma(a+1)\Gamma(b+1)$,  $\lim_{x\to 1^-}(1-x)\,_2F_1(a,b,a+b,x)^2=0,$ and
$$\lim_{x\to 1^-}\,_2F_1(a,b,a+b,x)-\frac{ab}{a+b}\,_2F_1(a,b,a+b+1,x)(c-\log(1-x))=\frac{R(a,b)-c}{B(a,b)},$$
therefore, a necessary condition for which $\varphi_{a,b,c}(x)$ is convex: $c\leq R(a,b)$.

3)
Let $U_{a,b}(x)=c-\log(1-x)$, $V_{a,b}(x)=\,_2{F}_1(a,b;a+b;x)$, $T_{a,b}(x)=U'_{a,b}(x)V_{a,b}(x)-U_{a,b}(x)V'_{a,b}(x)$, $Y_{a,b}(x)=V_{a,b}(x)^2$ and
$$Z_{a,b}(x)=\frac{T'_{a,b}(x)Y_{a,b}(x)}{Y_{a,b}'(x)}-T_{a,b}(x).$$
Then $$\varphi_{a,b}(x)=\frac{V_{a,b}(x)}{U_{a,b}(x)},\quad \varphi'_{a,b}(x)=\frac{T_{a,b}(x)}{Y_{a,b}(x)},$$
and \begin{equation}\label{e}\varphi''_{a,b}(x)=\frac{Z_{a,b}(x)Y_{a,b}'(x)}{Y_{a,b}(x)^2}.\end{equation}
Since, $$Y'_{a,b}(x)=2\frac{ab}{a+b}\,_2{F}_1(a+1,b+1;a+b+1;x)\,_2{F}_1(a,b;a+b;x)>0.$$
Differentiate yields
$$\frac{T'_{a,b}(x)}{Y_{a,b}'(x)}=\frac{B_{a,b}(x)}{A(x)}.$$
Then, $$Z'_{a,b}(x)=\Big(\frac{B_{a,b}(x)Y_{a,b}(x)}{A(x)}-T_{a,b}(x)\Big)'=\Big(\frac{B_{a,b}(x)}{A(x)}\Big)'Y_{a,b}(x)=\frac{A'(x)}{A(x)^2}P_{a,b}(x)Y_{a,b}(x)<0.$$
Therefore, and by using equation \eqref{e}, we get
$Z_{a,b}(x)<Z_{a,b}(0)=\varphi''_{a,b}(0)\frac{Y_{a,b}(0)^2}{Y_{a,b}'(0)}.$
Since, $Y_{a,b}'(0)>0$ and the condition on the parameter $c\geq c(a,b)$ is equivalent to  $\varphi''_{a,b}(0)<0$, therefore $Z_{a,b}(x)\leq Z_{a,b}(0)<0.$

Then $\varphi''_{a,b}(x)<0$. This completes the proof.

Next, assume that $\varphi_{a,b}$ is strictly concave on $(0,1)$, then $\varphi''_{a,b}(0)<0$. By straightforward computation, one gets
$$\varphi''_{a,b}(0)=\frac{a+b-2 a b}{a+b}-c\frac{a b(a+b+1)(a+b-2ab)+ab(a+b)}{(a+b)^2(a+b+1)}.$$
Then $$c\geq \frac{(a+b)(a+b+1)(a+b-2 a b)}{a b\Big((a+b+1)(a+b-2ab)+ab(a+b)\Big)}.$$

\subsection{Proof of Theorem \ref{t1}.}

$$f_{a,b,c}(x)=\frac{ _2F_1(a,b,a+b,x)}{c-\log(1-x)},$$
$$f_{a,b,c}'(x)=\frac{(ab/a+b) _2F_1(a,b,a+b+1,x)(c-\log(1-x))-\,_2F_1(a,b,a+b,x)}{(1-x)(c-\log(1-x))^2}=\frac{u(x)}{v(x)}$$

$$u'(x)=\frac{(ab)^2}{(a+b)(a+b+1)}\,_2F_1(a+1,b+1,a+b+2,x)(c-\log(1-x)),$$
$$v'(x)=-(c-\log(1-x))^2+2(c-\log(1-x)),$$
$$u'v=\frac{(ab)^2}{(a+b)(a+b+1)}\,_2F_1(a+1,b+1,a+b+2,x)(1-x)(c-\log(1-x))^2,$$
$$uv'=\Big((ab/a+b) _2F_1(a,b,a+b+1,x)(c-\log(1-x))-\,_2F_1(a,b,a+b,x)\Big)((c-\log(1-x))-2).$$
Then, $$f_{a,b,c}''(x)=\frac{\Psi_{a,b,c}(x)}{(1-x)^2(c-\log(1-x))^3},$$
where,
$$\Psi_{a,b,c}(x)=h(x)(c-\log(1-x))^2+g(x)(c-\log(1-x))+2\,_2F_1(a,b,a+b,x)$$

Using the notations of Lemma \ref{i}, we get $$\Psi(x)=h(x)(c-\varphi_+(x))(c-\varphi_-(x)).$$
Then, the function $f_{a,b,c}$ is convex on $(0,1)$ or equivalently $\Psi_{a,b,c}(x)\geq 0$ if and only if $c\geq\max_{[0,1)}\varphi_+(x)$ and is concave if and only if $\max_{x\in [0,1)}\varphi_-(x)\leq c\leq\min_{x\in[0,1)}\varphi_+(x)$

If there is $x_0\in(0,1)$ such that $\varphi_1(x_0)=\varphi_2(x_0)$, then $\Delta(x_0)=0$, since the function $x\mapsto\Delta(x)$ increases and nonnegative, then $0\leq \Delta(0)\leq\Delta(x_0)=0$, whence $a=b=1/2$ and $\varphi_1(0)=\varphi_2(0)=8/3$. In that case, the function is concave if and only if $c=8/3$

\section{Inequalities}
We make use of our convexity result to give
certain functional inequalities.
In [4], Theorem 1.1 and Theorem 1.2 was applied to obtain several functional inequalities involving the complete
elliptic integral $\cal K$ (see [4, Corollaries 4.1-4.4]) and for the generalized elliptic integral ${\cal K}_p$. It is worth noting that arguments
similar to those used there, together with the main results of this paper, can be applied to
obtain functional inequalities involving the Gauss hypergeometric function. For example, we present the following extensions of [4, Corollary 4.1]:
\begin{corollary}
For all $x\in(0,1)$, $a,b\in[0,1]$, $a+b\leq 1$ and $c\geq c(a,b)$
$$c+B(a,b)\leq \frac{c-\log(1-x)}{\,_2F_1(a,b,a+b,x)}+\frac{c-\log(x)}{\,_2F_1(a,b,a+b,1-x)}\leq\frac{(2c+2\log 2)}{\,_2{F}_1(a,b;a+b;1/2)}.$$
If $c\leq R(a,b)$ then the reversed inequalities holds.
\end{corollary}
\begin{corollary}
For all $x\in(0,1)$, $a,b\in[0,1]$, $a+b\leq 1$ and $c\in [\delta_-,\delta_+]$
$$\frac1c+\frac1{B(a,b)}\leq \frac{\,_2F_1(a,b,a+b,x)}{c-\log(1-x)}+\frac{\,_2F_1(a,b,a+b,1-x)}{c-\log(x)}\leq\frac{\,_2{F}_1(a,b;a+b;1/2)}{(2c+2\log 2)}.$$
If $c\geq \alpha_0$ then the reversed inequalities holds.
\end{corollary}
{\bf Proof.} Let $f(x)=G_{a,b}(x)+G_{a,b}(1-x)$, then
$$f'(x)=G'_{a,b}(x)-G'_{a,b}(1-x),\quad f''(x)=G''_{a,b}(x)+G''_{a,b}(1-x),$$
since, $G_{a,b}$ is concave the $f'(x)$ decreases, moreover, $f'(1/2)=0$. Hence, $f$ decreases on $(0,1/2)$ and increases on $(1/2,1)$. Then
$$2G_{a,b}(1/2)\leq G_{a,b}(x)<\max(G_{a,b}(0), G_{a,b}(1)),$$
Further, $G_{a,b}(1/2)=(c+\log 2)/\,_2{F}_1(a,b;a+b;1/2)$, $G_{a,b}(0)=c$ and $G_{a,b}(1)=\Gamma(a+b)/\Gamma(a)\Gamma(b)$. Thus,
$$c+\frac1{B(a,b)}\leq \frac{c-\log(1-x)}{\,_2F_1(a,b,a+b,x)}+\frac{c-\log(x)}{\,_2F_1(a,b,a+b,1-x)}\leq2(c+\log 2)/\,_2{F}_1(a,b;a+b;1/2)$$
\begin{proposition}\label{pro1} For $b>0$, and $x\in(0,1)$ let

$$F(x)=\frac{B(1/2+x,3/2-x)}{x(1-x)}+R(1/2+x,3/2-x)-R(x,1-x),$$
and

$$f_{b}(x)=(\frac1x+\frac1b)B(1/2+x,1/2+b)+R(1/2+x,1/2+b)-R(x,b).$$
\begin{enumerate}
\item The function $F$ is strictly completely monotonic $(0,1/2)$ and absolutely monotonic on $(1/2,1)$,
and $$4-4\log2<F(x).$$
Equality holds if and only if $x=1/2$
\item The function $f_{b}$ is strictly completely monotonic on $(0,1-b)$ for all $b\in(0,1)$.
\end{enumerate}
\end{proposition}
{\bf Proof.}
1) We have, $F(1-x)=F(x)$.
We set $$H(x)=\frac{\Gamma(3/2-x)\Gamma(1/2+x)}{x}+\psi(x)-\psi(\frac12+x),$$
then $$F(x)=H(x)+H(1-x).$$
We prove that $H$ is completely monotonic on $(0,1/2)$ and absolutely monotonic on $(1/2,1)$.
Using [26, Corollary 3.3], and the series expansion
\begin{equation}\label{ee2}\frac\pi{\sin(\pi z)}=\frac1z+2\sum_{n=1}^\infty\eta(2n)z^{2n-1},\quad|z|<1,\end{equation}
where $\eta(n)=\sum_{k=1}^\infty(-1)^{k-1}/k^n$ and $\eta(2n)=(1-2^{1-2n})\zeta(2n)$, here $\zeta$ states for The Riemann zeta function.
Since, $B(\frac12-x,\frac12+x)=\pi/\sin(\pi(1/2-x))$, therefore, for $x\in(0,1/2)$, we get
$$B(\frac12-x,\frac12+x)=\frac1{\frac12-x}+2\sum_{n=1}^\infty\eta(2n)(\frac12-x)^{2n-1},$$
and
$$\frac{\Gamma(3/2-x)\Gamma(1/2+x)}{x}=\frac1{x}\Big(\sum_{n=0}^\infty\widetilde\eta(n)(\frac12-x)^{n}\Big),$$
where $\widetilde\eta(0)=1$, $\widetilde\eta(1)=0$, $\widetilde\eta(2n)=2\,\eta(2n)$ and $\widetilde\eta(2n+1)=0$ for $(n\geq 1)$.
Since, $1/x=2\sum_{n=0}^\infty 2^n(1/2-x)^n$. Using the product of series for $x\in(0,1/2)$, we get
$$\frac{\Gamma(3/2-x)\Gamma(1/2+x)}{x}=\sum_{n=0}^\infty\alpha_n(\frac12-x)^{n},$$
where, $\alpha_{n}=\sum_{k=0}^n2^{n-k+1}\widetilde\eta(k)$.
Moreover, for $x>0$, $$\psi(x)=\int_0^\infty(\frac {e^{-u}}u-\frac{e^{-x u}}{1-e^{-u}})du.$$
Then, $$\psi(x)-\psi(\frac12+x)=-\int_0^\infty\frac {e^{-x u}}{1+e^{-\frac u2}}du=\sum_{n=0}^\infty\frac{a_n}{n!}(\frac12-x)^n,$$
where, $$a_n=-\int_0^\infty\frac{e^{-u/2}u^n}{1+e^{-\frac u2}}du.$$
Since, $$a_n=\sum_{k=1}^\infty(-1)^k\int_0^\infty e^{-ku/2}u^ndu=2^{n+1}n!\sum_{k=1}^\infty\frac{(-1)^k}{k^{n+1}}.$$
Then, $$\psi(x)-\psi(\frac12+x)=-2^{n+1}\sum_{n=0}^\infty\eta(n+1)(\frac12-x)^n,$$
and $$H(x)=\sum_{n=0}^\infty\gamma_n(\frac12-x)^n.$$
where $$\gamma_n=\alpha_n-2^{n+1}\eta(n+1),$$
then $$\gamma_n=\sum_{k=0}^n2^{n-k+1}\widetilde\eta(k)-2^{n+1}\sum_{k=1}^\infty\frac{(-1)^{k-1}}{k^{n+1}},$$
or $$\gamma_n=\sum_{k=0}^n2^{n-k+1}\widetilde\eta(k)+(\sum_{k=1}^\infty\frac{1}{k^{n+1}}-\sum_{k=1}^\infty\frac{1}{(k+\frac12)^{n+1}}).$$
Since, for all $n$, $\widetilde\eta(n)\geq 0$, then $\gamma_n>0$ for all $n$.

Therefore, for $x\in(0,1/2)$ and all $m$
$$(-1)^mH^{(m)}(x)=\sum_{n=m}^\infty\frac{(n-k)!}{n!}\gamma_n(\frac12-x)^{n-m}>0,$$
and $H$ is strictly completely monotonic on $(0,1/2)$.

On the other hand and by using the formulas: for $x>0$
$$\psi(1+x)=\psi(x)+\frac1x,\quad\psi(1-x)=\psi(x)+\pi\cot(\pi x).$$
It follows that
$$H(1-x)=H(x)+\frac{2\pi}{\sin(2\pi x)}+\frac1{x-\frac12},$$
or,
$$H(1-x)=H(x)+\frac{2\pi}{\sin(2\pi (\frac12-x))}-\frac1{\frac12-x}.$$
Moreover, for $x\in(0,1/2)$, and from equation \eqref{ee2}, we get
$$\frac{2\pi}{\sin(2\pi (\frac12-x))}-\frac1{\frac12-x}=\sum_{n=1}^\infty2^{2n+1}\eta(2n)(1/2-x)^{2n-1},$$
and for all $n$, $\eta(2n)>0$, hence, the function $x\mapsto \frac{2\pi}{\sin(2\pi (\frac12-x))}-\frac1{\frac12-x}$ is strictly completely monotonic on $(0,1/2)$. Using the fact that $H$ is strictly completely monotonic, we get $H^{(m)}(1-x)>0$ for all $x\in(0,1/2)$. This implies that $H$ is absolutely monotonic on $(1/2,1)$.

Now, $F(x)=H(x)+H(1-x)$, then for all $m$ and all $x\in(0,1/2)$
$$(-1)^mF^{(m)}(x)=(-1)^mH^{(m)}+H^{(m)}(1-x)>0.$$
From the equation $F(x)=F(1-x)$, we get for all $x\in(1/2,1)$ and $m\in\Bbb N$, $F^{(m)}(x)=(-1)^mF^{(m)}(1-x)>0$.

2)
Recall that the beta function admits the following integral rpresentation
$$B(x,y)=\frac{\Gamma(x)\Gamma(y)}{\Gamma(x+y)}=\int_0^\infty e^{-x u}(1-e^{-u})^{y-1}du.$$
$$f_1(x):=\frac{\Gamma(x+1/2)\Gamma(b+1/2)}{b x\Gamma(x+b)}=\int_0^\infty \theta(t) e^{-x t}dt,$$
where, $$\theta(t)=\frac1b(1-e^{-t})^{b-1/2}e^{-t/2}+\int_0^t(1-e^{-u})^{b-1/2}e^{-u/2} du.$$
Since,
$$f'_{b}(x)=f_1'(x)-\psi'(x+1/2)+\psi'(x)=\int_0^\infty\varphi_{b}(t)e^{-x t},$$
where, $$\varphi_{b}(t)=-t\theta(t)+\frac{t}{1+e^{-t/2}}.$$
Let $$g(t,b)=\frac1b(1-e^{-t})^{b-1/2}e^{-t/2}+\int_0^t(1-e^{-u})^{b-1/2}e^{-u/2} du,$$
then
$$\frac{\partial}{\partial b} g(t,b)=(-\frac1{b^2}+\log(1-e^{-t}))(1-e^{-t})^{b-1/2}e^{-t/2}+\int_0^t(1-e^{-u})^{b-1/2}e^{-u/2}\log(1-e^{-u}) du<0,$$
therefore, $b\mapsto g(t,b)$ is strictly decreasing. Then the function $b\mapsto  -\theta(t)+\frac{1}{1+e^{-t/2}}$ is strictly increasing and then
$\varphi_b(t)>\varphi_{1}(t).$ Moreover,
$$\varphi_1(t)=t(-(1-e^{-t})^{1/2}e^{-t/2}-\int_0^t(1-e^{-u})^{1/2}e^{-u/2} du+\frac{1}{1+e^{-t/2}}):=t\Theta(t).$$
Since,
$$\begin{aligned}\Theta'(t)e^{-t/2}&=-\frac12(1-e^{-t})^{1/2}-\frac12(1-e^{-t})^{-1/2}e^{-t}+\frac12\frac{1}{(1+e^{-t/2})^2}\\
&\leq -\frac12(1-e^{-t})^{1/2}-\frac12(1-e^{-t})^{-1/2}e^{-t}+\frac12\frac{(1-e^{-t/2})^{1/2}}{(1-e^{-t})^{1/2}}\\
&=\frac1{2(1-e^{-t})^{1/2}}(-1+(1-e^{-t/2})^{1/2})<0.
\end{aligned}$$
Therefore, $(-1)^nf^{(n)}_{b}(x)>0$ for all $n\geq 1$. We deduce that $f_b(x)$ is strictly decreasing and for all $x\in(0,1-b)$,
$f_b(x)>f_{b}(1-b)=F(b)>0$ by the first item. Whence the function $f_b(x)$ is strictly completely monotonic on $(0,1-b)$ for all $b\in(0,1)$.

\begin{proposition} For $a,b\geq 0$ and $a+b\leq 1$, let
$$G_{a,b,c}(x)=\frac{c+ x\,_2F_1(a+1/2,b+1/2,a+b+1,x)}{\,_2F_1(a,b,a+b,x)},$$
The function $H_{a,b,c}$ is strictly decreasing on $(0,1)$ if and only if $c\geq 1/a+1/b$, and is strictly increasing on $(0,1)$ if and only if
$$c<\frac{R(a,b)-R(a+1/2,b+1/2)}{B(1/2+a,1/2+b)}.$$
\end{proposition}
For $a=1/2, b=1/2$, we get $G_{1/2,1/2,c}(x)=(c-\log(1-x))/K(x)$, which corresponds to the result of Anderson et al, where $K$ is the elliptic function of
the first Kind.

{\bf Proof.} The function $G_{a,b,c}$ is symmetric on the variable $a,b$, then by assumptions one may suppose that $a\in(0,1/2)$ and $b\in(0,1-a)$.

1) Assume that $c\geq 1/a+ 1/b$ and let us write
$$G_{a,b,c}(x)=\frac{\sum_{n=0}^\infty a_nx^n}{\sum_{n=0}^\infty b_nx^n}=\frac{U(x)}{V(x)},$$
where $a_0=c$ and $a_n=(a+1/2)_{n-1}(b+1/2)_{n-1}/(n-1)!(a+b+1)_{n-1}$, and $b_n=(a)_{n}(b)_{n}/n!(a+b)_{n}$, hence
$$u_n=\frac{a_n}{b_n}=\frac{n(a+b)(a+1/2)_{n-1}(b+1/2)_{n-1}}{(a)_{n}(b)_{n}}.$$
and by the formula $\Gamma(x+n+1)=(x+n)\Gamma(x+n)$, we get
$$\frac{u_{n+1}}{u_{n}}=(1+\frac1n)(1-\frac1{2(a+n)})(1-\frac1{2(b+n)}).$$
The function $b\mapsto 1-1/(2(b+n))$ is strictly increasing, then
for $a\in(0,1/2)$ and $b\in(0,1-a)$
$$\frac{u_{n+1}}{u_{n}}\leq(1+\frac1n)(1-\frac1{2(a+n)})(1-\frac1{2(1-a+n)}):=v_n(a).$$
Differentiate with respect to $a$ yields
$$\begin{aligned}v'_n(a)&=\frac{-4 a^4+8 a^3+a^2 \left(12 n^2+8 n-5\right)+a \left(-12 n^2-8 n+1\right)+2 n (n+1)^2}{4 n^2 (-a+n+1)^2 (a+n)^2}\\
&:=\frac{w_n(a)}{4 n^2 (-a+n+1)^2 (a+n)^2}.\end{aligned}$$
Moreover, The derivatives of the function $w_n$ with respect to $a$ are
$w'(a)=(1-2 a)(8 a^2-8 a-12 n^2-8 n+1),\,w''(a)=48 a(1 - a) -10 +16 n + 24 n^2>0$ for all $n\geq 1$. So, $w'(a)\leq w'(1/2)=0$ and $w_n(a)$ is strictly decreasing on $(0,1/2)$ for all $n\geq 1$ and $w_n(a)>w_n(1/2)=(1 + 2 n)n^2$. This implies that that $v'_n(a)>0$ and $v_n(a)$ is strictly increasing. Then
$$\frac{u_{n+1}}{u_{n}}\leq(1+\frac1n)(1-\frac1{2n+1})^2.$$
By easy computation one shows that the sequence in the right-hand side is strictly increasing and $\displaystyle\lim_{n\to\infty} (1+\frac1n)(1-\frac1{2n+1})^2=1$. Then,
$u_{n+1}<u_n$ for all $n\geq 1$. By assumption $a_0/b_0=c>1/a+1/b=a_1/b_1$. Then the sequence $u_n=a_n/b_n$ decreases for all $n$. Applying Lemma \ref{lem2}, we deduce that the function $G_{a,b,c}$ is strictly decreasing.


Assume that $G_{a,b,c}$ is strictly decreasing, then $G'_{a,b,c}(x)<0$, therefore, $U'(0)V(0)-U(0)V'(0)=1-c (ab)/(a+b)<0,$ or else $c>1/a+1/b$.

2) Recall that for $n\geq 1$, $u_{n+1}<u_n$ and $u_1=1/a+1/b$ and $u_0=c$. Assume $a+b\leq 1$ and $$c<\frac{\Gamma(a+b+1)}{\Gamma(a+1/2)\Gamma(b+1/2)}(\psi(a+1/2)+\psi(b+1/2)-\psi(a)-\psi(b)),$$
 then by Proposition \ref{pro1}, we get $c<1/a+1/b$ then $u_0<u_1$.
 Moreover, $$U'(x)=\,_2F_1(a+1/2,b+1/2,a+b+1,x)+\frac{(a+1/2)(b+1/2)}{a+b+1}
x\,_2F_1(a+3/2,b+3/2,a+b+2,x),$$
$$V'(x)=\frac{ab}{a+b}\,_2F_1(a+1,b+1,a+b+1,x).$$
Then, $$H_{U,V}(x)=\frac{a+b}{ab}\frac{S(x)\,_2F_1(a,b,a+b,x)}{\,_2F_1(a,b,a+b+1,x)}
-c- x\,_2F_1(a+1/2,b+1/2,a+b+1,x),$$
where,
$$S(x)=(1-x)\,_2F_1(a+1/2,b+1/2,a+b+1,x)+\frac{(a+1/2)(b+1/2)}{a+b+1}
x\,_2F_1(a+1/2,b+1/2,a+b+2,x).$$
Using Lemma \ref{lem1}, we get for $x\to 1^-$
$$\begin{aligned}H_{U,V}(x)&=\frac1{B(a+1/2,b+1/2)}(R(a,b)-\log(1-x)+o((1-x)\log(1-x)))\\
&-c-\frac1{B(a+1/2,b+1/2)}(R(a+1/2,b+1/2)-\log(1-x)+o((1-x)\log(1-x)),\end{aligned}$$
Then, $$\lim_{x\to 1^-}H_{U,V}(x)=-c+\frac{\Gamma(a+b+1)}{\Gamma(a+1/2)\Gamma(b+1/2)}(R(a,b)-R(a+1/2,b+1/2)).$$

Since, $c> 1/a+1/b$, then by Proposition \ref{pro1} we get $\lim_{x\to 1^-}H_{U,U}(x)< 0$, and
 From Lemma \ref{lem2}, we deduce that the function $U(x)/V(x)$ is strictly increasing.

Next, assume that $G_{a,b,c}$ is strictly increasing. We have,
$$G'_{a,b,c}(x)=\frac{V'(x)}{V(x)^2}H_{U,V}(x),$$
and $V'(x)>0$ for all $x\in(0,1)$, then $H_{U,V}(x)\geq 0$. By the previous computation, we get
$$-c+\frac{\Gamma(a+b+1)}{\Gamma(a+1/2)\Gamma(b+1/2)}(R(a,b)-R(a+1/2,b+1/2))\geq 0.$$
This completes the proof of the Proposition.

\begin{proposition} For $x\in(0,1)$, let
$$B(x)=B(x,1-x),\quad R(x)=R(x,1-x).$$
1) The function $f(x)=(1+x(1-x)-(x(1-x))^2)B(x)-R(x),$
is strictly completely monotonic on $(0,1/2)$.\\
2) The function $g(x)=R(x)-B(x)/(1+x(1-x))$ is strictly completely monotonic on $(0,1/2)$.
\end{proposition}
{\bf Proof.} We have
\begin{equation}\label{a}R(x)=\log 16+4\sum_{n=1}^\infty\lambda(2n+1)(1-2x)^{2n},\end{equation}
and \begin{equation}\label{b}B(x)=4\sum_{n=0}^\infty\beta(2n+1)(1-2x)^{2n},\end{equation}
where $$\lambda(n)=\sum_{k=0}^\infty\frac1{(2k+1)^n},\quad \beta(n)=\sum_{k=0}^\infty\frac{(-1)^k}{(2k+1)^n}.$$
Since, $1+x(1-x)-(x(1-x))^2=19/16-1/8(1-2x)^2-1/16(1-2x)^4,$ then
$$(1+x(1-x)-(x(1-x))^2)B(x)=\frac{19}4\beta(1)+(\frac{19}4\beta(3)-\frac12\beta(1))(1-2x)^2+\sum_{n=2}^\infty a_n(1-2x)^{2n},$$
where $$a_n=\frac{19}4\beta(2n+1)-\frac12\beta(2n-1)-\frac14\beta(2n-3).$$
$$f(x)=\sum_{n=0}^\infty b_n(1-2x)^{2n},$$
where $b_0=\frac{19}4\beta(1)-\log 16$, $b_1=\frac{19}4\beta(3)-\frac12\beta(1)-4\lambda(3)$ and for $n\geq 2$, $$b_n=\frac{19}4\beta(2n+1)-\frac12\beta(2n-1)-\frac14\beta(2n-3)-4\lambda(2n+1).$$
$b_0\simeq0.95$, $b_1\simeq 0.919$.
Since, $$b_n=\sum_{k=0}^\infty\frac{\alpha_k}{(2k+1)^{2n-3}},$$
with $$\alpha_k=\frac{(19/4)(-1)^k-4}{(2k+1)^{4}}-\frac12 \frac{(-1)^k}{(2k+1)^{2}}-\frac{(-1)^k}4.$$
Besides,
$$\begin{aligned}&\frac{\alpha_{2k}}{(4k+1)^{2n-3}}+\frac{\alpha_{2k-1}}{(4k-1)^{2n-3}}\\
&=\frac{\frac{3}{4(4k+1)^{4}}-\frac{1}{2(4k+1)^{2}}-\frac14}{(4k+1)^{2n-3}}+\frac{\frac{-35}{4(2k+1)^{4}}+\frac{1}{2(4k-1)^{2}}+\frac{1}4}{(4k-1)^{2n-3}}\\
&=\frac{(4k-1)^4+2(4k-1)^2-35}{4(4k-1)^{2n+1}}-\frac{(4k+1)^4+2(4k+1)^2-3}{4(4k+1)^{2n+1}}\\
&=\frac{((4k-1)^2-5)((4k-1)^2+7)}{4(4k-1)^{2n+1}}-\frac{((4k+1)^2-1)((4k+1)^2+3)}{4(4k+1)^{2n+1}}.\end{aligned}$$
For $x>0$ and $y>0$ let
$$\varphi_x(y)=\frac{(x^2-5)(x^2+7)}{4x^{2y+1}}-\frac{(x^2+4x+3)(x^2+4x+7)}{4(x+2)^{2y+1}},$$
then
$$\varphi_x'(y)=2 (x^2+4 x+3)(x^2+4x+7) (x+2)^{-2 y-1} \log (x+2)-2(x^2-5)(x^2+7) x^{-2 y-1} \log (x).$$
So, for $x\geq 5$,
$$\varphi_x'(y)\leq 2 \theta(x) x^{-2 y-1}\log (x),$$
where $$\theta(x)=(x^2+4 x+3)(x^2+4x+7)\log (\frac75)-(x^2-5)(x^2+7).$$
Differentiate yields
$$\theta'(x)=4 ((\sigma-1)x^3+6 \sigma x^2+(12 \sigma-1)x+8\sigma),$$
and
$$\theta''(x)=4(3(\sigma-1)x^2+12 \sigma x +12 \sigma-1),$$
$\sigma=\log(7/5)$. For $x\geq 5$, $\theta''(x)<0$ and
$\theta'(5)<0$, then $\theta$ is strictly decreasing and $\theta(7)>0$. Thus $\varphi'_x(y)>0$ for $x\geq 7$. Therefore, $\varphi_x$ is strictly decreasing and $\lim_{y\to +\infty}\varphi_x(y)=0$ for all $x\geq 7$. Therefore,
$$\frac{\alpha_{2k}}{(4k+1)^{2n-3}}+\frac{\alpha_{2k-1}}{(4k-1)^{2n-3}}=\varphi_{4k-1}(n)>0,$$
for all $k\geq 2$, and $n\geq 1$.

Furthermore, $\alpha_0=0$ and $\alpha_1+\alpha_2=8(2/3^{2n+1}-21/5^{2n+1})>0$ for all $n\geq 2$, and
$$b_n=\sum_{k=1}^\infty(\frac{\alpha_{2k-1}}{(4k-1)^{2n-3}}+\frac{\alpha_{2k}}{(4k+1)^{2n-3}})>0$$ for all $n\geq 2$ and $b_0\simeq0.95$, $b_1\simeq 0.919$. Therefore, the function $f$ is strictly completely monotonic.

2) Firstly $1+x(1-x)=\frac54-\frac14(1-2x)^2),$ then for $x\in(0,1/2)$
$$\frac1{1+x(1-x)}=\frac45\sum_{n=0}^\infty\frac1{5^n}(1-2x)^{2n},$$
From equation \eqref{b}, we get
$$\frac{B(x)}{1+x(1-x)}=\frac{16}5\sum_{n=0}^\infty c_n(1-x)^{2n},$$
where $\displaystyle c_n=\sum_{m=0}^{n}\frac{\beta(2m+1)}{5^{n-m}}$. Therefore, and by using equation \eqref{a} we obtain
$$g(x)=\sum_{n=0}^\infty d_n(1-2x)^{2n},$$
where $d_0=\log16-(16/5)\beta(1)$ and $d_n=4(\lambda(2n+1)-(4/5) c_n)$ for $n\geq 1$.
$d_0=(4/5)\pi\log 16> 0.25$, and for $n\geq 1$
$$\frac14 d_n=\sum_{k=0}^\infty\Big(\frac1{(2k+1)^{2n+1}}-\frac45\sum_{m=0}^{n}\frac{(-1)^k}{5^{n-m}(2k+1)^{2m+1}}\Big).$$
Let
$$u_k(n)=\frac1{(2k+1)^{2n+1}}-\frac45\sum_{m=0}^{n}\frac{(-1)^k}{5^{n-m}(2k+1)^{2m+1}},$$
then for all $k\geq 1$
$$u_{2k-1}(n)+u_{2k}(n)=\frac1{(4k-1)^{2n+1}}+\frac1{(4k+1)^{2n+1}}+\frac45\sum_{m=0}^{n}\frac{1}{5^{n-m}}\frac{(4k+1)^{2m+1}-(4k-1)^{2m+1}}{(16k^2-1)^{2m+1}}>0,$$
and $u_0(n)=1-\frac45\sum_{m=0}^{n}\frac{1}{5^{n-m}}=1/5^{n+1}>0$.
Therefore, for all $n\geq 1$, $$\frac14 d_n=u_0(n)+\sum_{k=1}^\infty (u_{2k-1}(n)+u_{2k}(n))>0.$$

So, for all $m\in\Bbb N$ and $x\in(0,1/2)$
$$g^{(2m)}(x)=\sum_{n=2m}^\infty \frac{(2n)!}{(2n-2m)!}d_n(1-2x)^{2n-2m}>0,$$
and
$$-g^{(2m+1)}(x)=\sum_{n=m+1}^\infty \frac{(2n)!}{(2n-2m-1)!}d_n(1-2x)^{2n-2m-1}>0,$$
This implies that the function $g$ is strictly completely monotonic on $(0,1/2)$.
\begin{corollary} For all $x\in(0,1)$, we have
$$\log16-\frac{4\pi}5+\frac{B(x)}{1+x(1-x)}<R(x)<1+\frac{B(x)}{1+x(1-x)},$$
$$\log16 -\frac{4\pi} 5+(14\zeta(3)-\frac{2\pi}{25} (8 + 5 \pi^2))(x - 1/2)^2+\frac{B(x)}{1+x(1-x)}<R(x)<1+\frac{B(x)}{1+x(1-x)}.$$
$$(1+x(1-x)-(x(1-x))^2)B(x)-1<R(x)<(1+x(1-x)-(x(1-x))^2)B(x)-\frac{21\pi}{16}+\log16$$
\end{corollary}
\begin{problem} Find conditions on the parameters $a,b,c,d$ for which the function
$$G_{a,b,c,d,\alpha}(x)=\frac{\alpha+ x\,_2F_1(c,d,c+d,x)}{\,_2F_1(a,b,a+b,x)},$$
is increasing (resp. decreasing) on $(0,1)$ .
\end{problem}
\bibliographystyle{plain}
\bibliography{ellipticfunction}

\end{document}